\documentclass[10pt,a4paper,titlepage]{article}
\usepackage[latin1]{inputenc}
\usepackage{amsmath}
\usepackage[ps,all,arc,rotate]{xy}
\usepackage{amsfonts}
\usepackage{amssymb}
\usepackage{amsthm}
\usepackage[dvips]{graphicx}
\usepackage{fancyhdr}
\usepackage{float}
\author{Dawid Kielak}
\title{Groups with infinitely many ends are not fraction groups}
\makeatletter
\def\imod#1{\allowbreak\mkern10mu({\operator@font mod}\,\,#1)}
\makeatother
\numberwithin{figure}{section}
\theoremstyle{plain}
\newtheorem{thm}{Theorem}[section]
\newtheorem*{prop*}{Proposition}
\newtheorem*{thm*}{Theorem}

\newtheorem{lem}[thm]{Lemma}
\newtheorem{cor}[thm]{Corollary}
\theoremstyle{definition}

\newtheorem{dfn}[thm]{Definition}
\newtheorem*{dfn*}{Definition}

\theoremstyle{remark}

\def\iff{if and only if }

\newdir{ >}{{}*!/-5pt/@{>}}
\newdir{< }{{}*!/+5pt/@{<}}

\begin{document}
\textsc{\begin{LARGE}\begin{center} Groups with infinitely many ends are not fraction groups
\end{center}\end{LARGE}}

\medskip

\begin{center}
Dawid Kielak\footnotemark

University of Bonn

\today
\end{center}

\medskip

\begin{center}
\begin{minipage}{0.65\textwidth}
\textsc{Abstract.} We show that any finitely generated group $F$ with infinitely many ends is not a group of fractions of any finitely generated proper subsemigroup $P$, that is $F$ cannot be expressed as a product $P P^{-1}$. In particular this solves a conjecture of Navas in the positive.
As a corollary we obtain a new proof of the fact that finitely generated free groups do not admit isolated left-invariant orderings.
%We also exhibit a new class of groups admitting an isolated left-invariant ordering, whose positive cone is not finitely generated.
\end{minipage}
\end{center}

\bigskip

\footnotetext[1]{Author supported by the ERC Grant Nb. 10160104}

\section{Introduction}

The existence of a left-invariant order on a group $G$ is equivalent to the existence of a positive cone $P \subset G$, that is a subsemigroup such that $G$ can be written as a disjoint union $G = \{1\} \sqcup P \sqcup P^{-1}$. In fact there is a one-to-one correspondence between left-invariant orderings and such positive cones.

In this note we prove that whenever a finitely generated group $F$ with infinitely many ends can be written as $F = P P^{-1}$, where $P$ is a finitely generated subsemigroup of $F$, then $P = F$.
Our result answers a question of Navas, who conjectured that finitely generated free groups are not groups of fractions of finitely generated subsemigroups $P$ with $P \cap P^{-1} = \emptyset$.

As an application we obtain a new proof of the fact that the space of left-invariant orderings of a finitely generated free group (endowed with the Chabauty topology) does not have isolated points. This result follows from the work of McCleary~\cite{mccleary1985},
but appears in this form for the first time in the work of Navas~\cite{navas2010}. It is worth noting that our proof is the first geometric one.

We also deduce that the left-orderings of finitely generated groups with infinitely many ends do not have finitely generated positive cones. This was already known for free products of left-orderable groups by the work of Rivas~\cite{rivas2012}.

Our theorem complements a folklore result stating that whenever $\mathcal S$ is a finite generating set for a group $G$, and $G$ does not contain a free subsemigroup, then $G$ is a group of fractions of $P$, the semigroup generated by $\mathcal S$.

We should note here that finitely generated groups with infinitely many ends have been classified by Stallings~\cite{stallings1968, stallings1971}. They are precisely those fundamental groups of non-trivial graphs of groups with exactly one edge and a finite edge group, which are finitely generated and not virtually cyclic.

\medskip

\textbf{Acknowledgments:} the author wishes to thank Andr\'{e}s Navas for introducing him to this problem, and for his many comments. He also wishes to thank Thomas Haettel and Ursula Hamenst\"{a}dt for many helpful conversations, and the referee for pointing out ways of significantly improving the presentation of this note.

\section{The result}
In the following we will use $X$ to denote the (right) Cayley graph of a finitely generated group $F$ with respect to some finite generating set. We will identify $F$ with vertices of $X$, and use $d$ to denote the standard metric on the Cayley graph $X$. The isometric left-action of $F$ on $X$ and its subsets will be denoted by left multiplication. The notation $B(x,\xi)$ will stand for the closed ball centred at $x$ of radius $\xi$.

We will assume that $F$ has infinitely many ends, and so there will exist a constant $\kappa$ such that the ball $B = B(1,\kappa)$ disconnects $X$ into a space with at least 3 infinite components. We will use $S$ to denote the set of vertices of $B$.

\begin{dfn}
We say that $A \subset X$ is a \emph{shoot} \iff there exists $w \in F$ such that $A$ is a connected component of $X \smallsetminus w B$. We say that $w B$ \emph{bounds} the shoot.
%
%If $A$ and $A'$ are both shoots, and $A' \subseteq A$, then we say that $A'$ is a \emph{subshoot} of $A$.
\end{dfn}

\begin{lem}
\label{lem: placing ball}
Let $A$ be an infinite shoot bounded by $B$. Then there exists $w \in F$ such that $w(X \smallsetminus A) \subseteq A$ and $w^{-1}(X \smallsetminus A) \subseteq A$.
\begin{proof}
%Take $w_0 \in F$ such that $w_0 z \in A$ and $w_0 B \subset A$. If $w_0 B$ separates $w_0 z$ from the ball bounding $A$ then we are done. Suppose it is not the case. Take $A'$ to be a shoot bounded by $w_0 B$, which does not contain $w_0 z$. Note that $A' \subseteq A$.
%
Note that the ball $B(1,2 \kappa)$ is finite, since $X$ is locally finite.
Since $F$ has infinitely many ends and $A$ is infinite, there exists $\lambda$ such that
\[ L = \{ x \in F \mid d(1,x) = \lambda \} \cap A\]
 has more than $\vert B(1,2 \kappa) \vert$ elements.
Take $l \in L$. The cardinality of $L$ guarantees that there exists $l' \in L$ such that $l' l^{-1} \not\in B(1,2 \kappa)$.

Let $w = l' l^{-1}$. Observe that
\[d(w, l') = d(l' l^{-1},l')= d(1, l) = \lambda \]
Consider a shortest path between $w$ and $l'$. If it lies entirely in $A$, then in particular so does $w$. If not, then it must contain some point $b \in B$, since $B$ bounds $A$. Now we have
\[ \lambda = d(w,l') = d(w,b) + d(b,l') \geqslant d(w,b) + \lambda - \kappa\]
which implies that $d(w,b) \leqslant \kappa$, and hence that $w \in B(1,2 \kappa)$, which is a contradiction.
We have thus established that $w \in A \smallsetminus B(1,2 \kappa)$, and therefore that $wB \subset A$.

Note that $w^{-1} = l l'^{-1} \notin B(1,2 \kappa)$ enjoys the same properties as $w$, and so we immediately conclude that $w^{-1}B \subset A$, or equivalently that $B \subset wA$.

%Since $l \in A$, we have $1 \in l^{-1}A$ and thus we see that $l' \in l' l^{-1} A=wA$. Suppose that $1 \not\in wA$. Take a path between $1$ and $l'$ of length $\lambda$. Since $1 \not\in wA$, the path has to cross $wB$. But then
%\[ \lambda \geqslant d(wB, l') + d(wB, 1) \geqslant \lambda-\kappa + \kappa + 1 = \lambda+1,\]
%which is a contradiction. We conclude that $1 \in wA$. Again, since $w \not\in B(1,2 \kappa)$, we have $B \subset wA$.
%
%As $z \not\in A$, we have $wz \not\in wA$. But $z \in F \smallsetminus S$ and so $z$ lies in some shoot bounded by $B$. Thus also $wz$ lies in some shoot bounded by $wB$.
Since $wB$ and $B$ are disjoint, every shoot bounded by $wB$ either contains $B$ or is disjoint from it. Clearly, there is a unique shoot bounded by $wB$ containing $B$, and we have already shown that it is $wA$. Each of the other shoots bounded by $wB$ lies in a single shoot bounded by $B$, namely in the shoot bounded by $B$ which contains $wB$. But we have already seen that this is $A$. We are left with the conclusion that
$w(X \smallsetminus A) \subset A$, and our proof is finished by making the analogous observations for $w^{-1}$.
\end{proof}
\end{lem}

%Our lemma implies the following statement (which is by no means new, and follows directly from the work of Stallings).
%
%\begin{cor}
%Every finitely generated group with infinitely many ends contains a free subgroup of rank 2.
%\begin{proof}[Sketch proof]
%Let $A_0, A_1$ and $A_2$ denote three distinct infinite shoots bounded by $B$. Let $w_i$ for $i \in \{1,2\}$ be an element of $F$ with $w_i(X \smallsetminus A_i) \subset A_i$ and $w_i^{-1}(X \smallsetminus A_i) \subset A_i$, given by the lemma above. Now the Ping-Pong lemma applied to the action of $w_1$ and $w_2$ on $A_0$ shows that they generate a free subgroup.
%\end{proof}
%\end{cor}
%
We are now ready for the main result.

\begin{figure}
\centering
\includegraphics[scale=0.4]{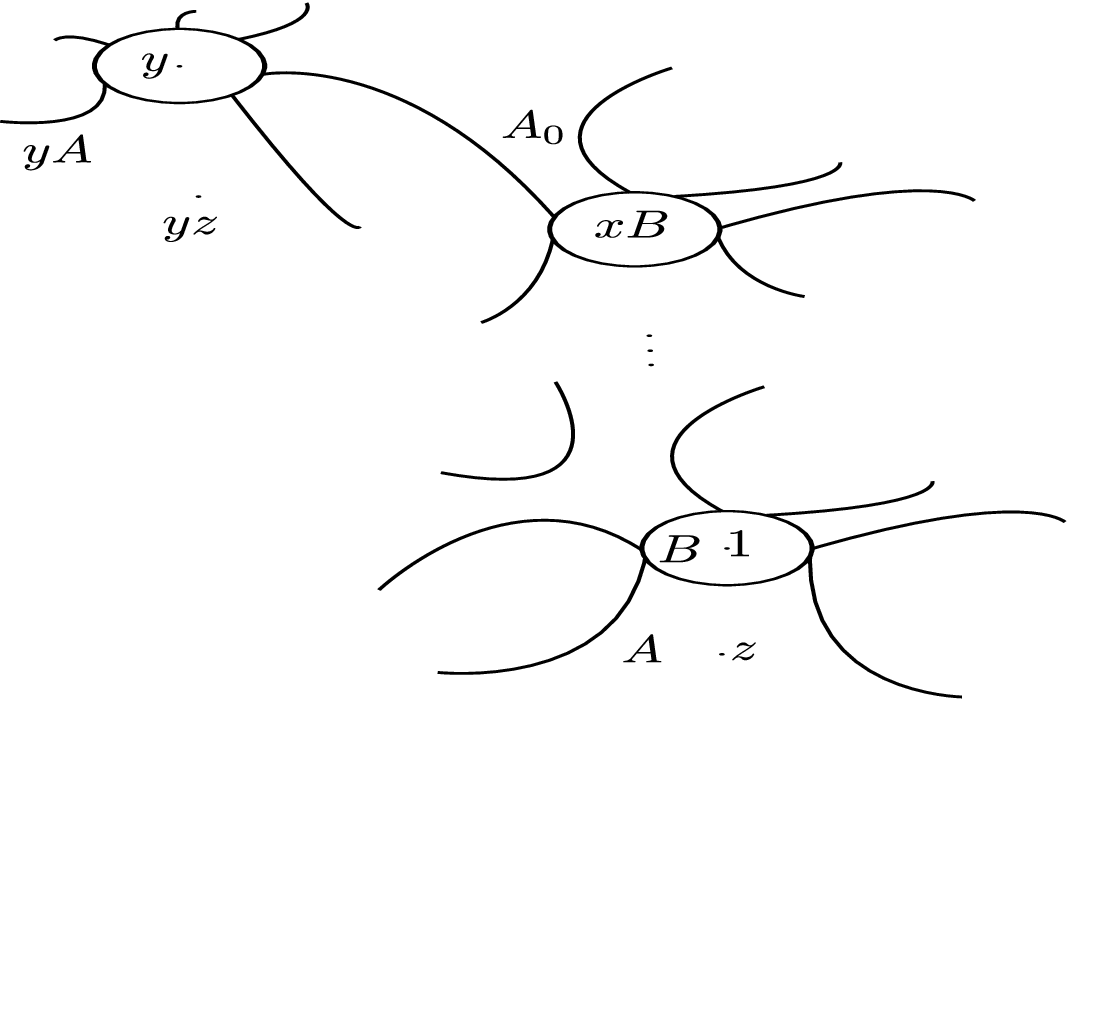}
\vspace{-20mm}
\caption{Step 1 of the theorem.}
\label{fig: figure}
\end{figure}

\begin{thm}
\label{thm: result}
Let $P$ be a finitely generated subsemigroup of a finitely generated group $F$ with infinitely many ends. If
$P P^{-1} = F$
then $P=F$.
\begin{proof}
For ease of notation we will refer to the elements of $P$ as positive, and to the elements of $P^{-1} $ as negative.

We first note that any finite generating set of $P$ is a generating set for $F$.
Let $X$ be the Cayley graph of $F$ with respect to some such generating set.
Note that this allows us to view generators of $P$ as positive edges of $X$, and hence any positive element $p \in P$ is realised by a positive path between 1 and the vertex $p$ in $X$.

We will use the notation $\kappa$, $B$ and $S$ as defined above.

\medskip
\noindent \textbf{Step 1:} We claim that $S (P^{-1} \cup \{1\}) = F$.
\medskip

If $P$ intersects each ball $B(x,\kappa)=xB$ then each $x \in F$ is a concatenation of an element of $P$ (namely any positive path from $1$ to $xB$) with an element in $S$ (connecting the end of the positive path to the centre of the ball). Thus we have $x \in P S$, and our claim follows by taking inverses.

Let us now suppose that there exists an $x \in F$ such that
\[P \cap xB = \emptyset\]
Let $A_0$ denote an infinite shoot bounded by $xB$ such that $1 \not\in A_0$.

Let $z \in F \smallsetminus S$ be any element, and let $A$ be the shoot bounded by $B$ containing $z$. We claim that there exists $y \in F$ such that $yA \subseteq A_0$.

There are two cases we need to consider. The first one occurs when
\[xA \subseteq A_0\]
in which case we take $y=x$. The other one (illustarted in Figure~\ref{fig: figure}) occurs when $xA \not\subseteq A_0$, that is when $xA$ is a shoot bounded by $xB$ other than $A_0$.
Lemma~\ref{lem: placing ball} applied to $x^{-1} A_0$ gives us an element $w \in F$ such that $w (X \smallsetminus x^{-1} A_0) \subseteq x^{-1} A_0$. So $y = x w$ satisfies
\[ yA = xwA \subseteq x x^{-1} A_0 = A_0 \]
and so we have proven the claim.

Now, since $yz \in F = P P^{-1}$, we can write $yz=pq$, where $p$ is positive and $q$ is negative. Since there are no positive elements in $xB$ by assumption, we see that $p \not\in A_0$, and therefore $q$ is a negative path connecting a vertex $p \in X \smallsetminus A_0$ to $yz\in yA \subseteq A_0$.
The shoot $yA$ is bounded by $yB$ and contained in $A_0$, hence any path from $X \smallsetminus A_0$ to $yA$ has to cross $yB$. This is in particular true for $q$,
so there is a negative path (a terminal subpath of $q$) from some vertex of $yB$ to $yz$, and hence from a vertex of $B$ to $z$ (after translating by $y^{-1}$). In the group language we have thus shown that $z \in S P^{-1}$, and so
\[ F \smallsetminus S \subseteq S P^{-1}\]
But clearly $S \subset S (P^{-1} \cup \{1\})$, and so we have proven the claim of step 1.

\medskip
\noindent \textbf{Step 2:} We claim that $P= F$.
\medskip

We have established above that $S (P^{-1} \cup \{1\})=F$, with $S$ being finite. Let $Q$ be a minimal (with respect to cardinality) finite subset of $F$ such that $Q (P^{-1} \cup \{1\})=F$. Suppose that there exist distinct $q,q' \in Q$. Then $q^{-1} q' \in F = P P^{-1}$, and so $q^{-1} q' = a b^{-1}$ with $a,b \in P$. Hence
\[q,q' \in qa P^{-1}\]
and therefore we could replace $Q$ by $(Q \cup \{qa\}) \smallsetminus \{q,q'\}$ of smaller cardinality. This shows that $\vert Q \vert = 1$. Without loss of generality we can take $Q = \{ 1 \}$, and thence get
\[P^{-1} \cup \{1\}=F\]
Now let $f \in F \smallsetminus \{1\}$. We have $f,f^{-1} \in P^{-1}$, and since $P^{-1}$ is a semigroup, also $1 = f f^{-1} \in P^{-1}$. So $P^{-1} = F$. Taking an inverse concludes the theorem.
\end{proof}
\end{thm}

We now easily deduce the following.

\begin{cor}
Let $F$ be a finitely generated group with infinitely many ends. Then $F$ does not allow a left-invariant ordering with a finitely generated positive cone.
\begin{proof}
Let $P$ be the positive cone of a left-invariant ordering of $F$. Then
\[F = P \cup P^{-1} \cup \{1\}\]
and so in particular $F=P P^{-1}$. But also $P \cap P^{-1} = \emptyset$, and so $P \neq F$. Now the contrapositive of Theorem~\ref{thm: result} tells us that $P$ is not finitely generated.
\end{proof}
\end{cor}

The statement above follows from
the work of Rivas~\cite{rivas2012}, since left-orderable groups are torsion-free, and
so they have infinitely many ends only when they are free products.

\medskip

%On the other hand it gives a new class of examples of groups with isolated orders whose positive cones are not finitely generated -- namely the class of finitely generated groups with infinitely many ends which admit an isolated order.
%It is however not clear, whether this class is non-empty.

We also get the following corollary.

\begin{cor}
The space of left-invariant orderings on any finitely generated free group has no isolated points.
\begin{proof}
Let $P$ be the positive cone of an isolated ordering of $F$, a finitely generated free group. By above, $P$ is not finitely generated.

The order defined by $P$ is isolated, and so there exists a finite set $S \subset F$ such that whenever we have another positive cone of an ordering $P'$ such that $P \cap S = P' \cap S$, then $P=P'$. However the work of Smith and Clay \cite[Theorem E]{claysmith2009} allows us to construct an order (in fact infinitely many such orders) whose positive cone $P'$ satisfies $P \cap S = P' \cap S$, but such that $P \neq P'$. This is a contradiction.
\end{proof}
\end{cor}

\medskip \textbf{Added in proof.} From the main theorem one can also easily deduce that groups with finite Garside structures have at most two ends.

\bibliographystyle{abbrv}
\bibliography{bibliography}

\noindent Dawid Kielak \newline
Mathematisches Institut der Universit\"at Bonn \newline
Endenicher Allee 60 \newline
D-53115 Bonn \newline
Germany \newline
\texttt{kielak@math.uni-bonn.de}

\end{document}